\title{A new construction relating enriched categories and internal ones in an extensive ambient}
\author{Matteo Di Domenico
        }
\date{}
\documentclass{article}
\usepackage[left=3cm, right=1.5cm]{geometry}
\usepackage[italian, english]{babel}
\usepackage[T1]{fontenc}
\usepackage{amsmath, amsfonts,amsthm}
\usepackage{euler}

\usepackage{mathabx}

\usepackage{multicol}
\usepackage{microtype} %it  improves line filling with: font expansion:  it horizontally expands the characters in order to optimally fill each line; character protrusion: it lets some characters protrude into the margins (typically the hyphens and punctuation signs)
\usepackage[bottom]{footmisc}%footnotes options

\usepackage{xcolor}
\usepackage{hyperref}
\hypersetup{hidelinks, 
	colorlinks=true, citecolor=blue, anchorcolor=red, urlcolor=black}

\usepackage[shortlabels]{enumitem}
\usepackage[normalem]{ulem}% per testo barrato, comando: \sout{}

\usepackage[mathscr]{urwchancal}
\DeclareFontFamily{OT1}{pzc}{}
\DeclareFontShape{OT1}{pzc}{m}{it}{<-> s * [1.10] pzcmi7t}{}
\DeclareMathAlphabet{\mathpzc}{OT1}{pzc}{m}{it}

\usepackage{version}%Defines macros \includeversion{NAME} and \excludeversion{NAME}, each of which defines an environment NAME whose text is to be included or excluded from compilation
\excludeversion{comment}

\usepackage[autostyle,italian=guillemets]{csquotes}
\usepackage[bibstyle=numeric,citestyle=numeric,hyperref,backend=biber, sorting=ynt]{biblatex}
\usepackage{tikz-cd}

\usetikzlibrary{positioning}
\usetikzlibrary{calc}
\usetikzlibrary{external}
\usetikzlibrary{arrows.meta}
\tikzset{
	labl/.style={anchor=south, rotate=55, inner sep=.5mm}
}

\usepackage{tabto}

\makeatletter
\let\orgdescriptionlabel\descriptionlabel
\renewcommand*{\descriptionlabel}[1]{%
	\let\orglabel\label
	\let\label\@gobble
	\phantomsection
	\edef\@currentlabel{#1\unskip}%
	\let\label\orglabel
	\orgdescriptionlabel{#1}%
}
\makeatother %to insert with \ref the customised name of the item in description ambient

\usepackage{bm}

\usepackage{fancyhdr}

\def\cat#1{\ensuremath{{\bm{\mathsf{#1}}}}}
\def\ct#1{\ensuremath{\mathpzc{#1}}}
\def\dfn#1{{\bfseries\itshape#1}}
\def\pullback{\ar[dr,phantom,"\lrcorner", shorten=2mm,start anchor={[xshift=-2ex, yshift=2ex]north west}]}
\def\pullbackk{\ar[dr,phantom,"\lrcorner", shorten=4mm,start anchor={[xshift=1.5mm, yshift=2mm]north west}]}

\def\endo#1{\ensuremath{\mathsf{End}(#1)}}
\def\idem#1{\ensuremath{\mathsf{Idem}(#1)}}

\def\split#1{\ensuremath{\ct{Split}(\ct{#1})}}

\def\en{\ensuremath{\mathrm{\mathsf{En}}}}
\def\inter{\ensuremath{\mathrm{\mathsf{S_i}}}}
\def\vcat{\ensuremath{\bm{\ct{V}}\text{-}\cat{Cat}}}
\def\vcats{\ensuremath{\bm{\ct{V}}\text{-}\cat{Cat^s}}}
\def\vcatk{\ensuremath{\bm{\ct{V}}\text{-}\cat{Cat^\kappa}}}
\def\catv{\ensuremath{\cat{Cat}\bm{(\ct{V})}}}
\def\catvs{\ensuremath{\cat{Cat^s}\bm{(\ct{V})}}}
\def\catvk{\ensuremath{\cat{Cat^\kappa}\bm{(\ct{V})}}}
\newcommand{\catspv}{\cat{Cat_{sp}}\bm{(\ct{V})}}
\newcommand{\catspvs}{\cat{Cat_{sp}^s}\bm{(\ct{V})}}
\newcommand{\qu}[1]{``#1''}
\def\split#1{\ensuremath{\mathsf{Split}(\ct{#1})}}
\def\balpha{\ensuremath{\bm{\alpha}}}
\def\cbalpha{\ensuremath{\hat{\bm{\alpha}}}}

\newcommand{\natto}{%
	\mathrel{\vbox{\offinterlineskip
			\mathsurround=0pt
			\ialign{\hfil##\hfil\cr
				\normalfont\scalebox{1.2}{.}\cr
				\noalign{\kern.05ex}
				$\to$\cr}
	}}%
}

\usepackage{pict2e}

\makeatletter
\newcommand{\adjunction}[4]{%
	#1\colon #2%
	\mathrel{\vcenter{%
			\offinterlineskip\m@th
			\ialign{%
				\hfil$##$\hfil\cr
				\longrightharpoonup\cr
				\noalign{\kern-.3ex}
				\smallbot\cr
				\longleftharpoondown\cr
			}%
	}}%
	#3 \noloc #4%
}
\newcommand{\longrightharpoonup}{\relbar\joinrel\rightharpoonup}
\newcommand{\longleftharpoondown}{\leftharpoondown\joinrel\relbar}
\newcommand\noloc{%
	\nobreak
	\mspace{6mu plus 1mu}
	{:}
	\nonscript\mkern-\thinmuskip
	\mathpunct{}
	\mspace{2mu}
}
\newcommand{\smallbot}{%
	\begingroup\setlength\unitlength{.15em}%
	\begin{picture}(1,1)
		\roundcap
		\polyline(0,0)(1,0)
		\polyline(0.5,0)(0.5,1)
	\end{picture}%
	\endgroup
}
\makeatother

\usepackage{graphicx}
\newtheoremstyle{plainnew}% name
{9pt}%      Space above, empty = `usual value'
{9pt}%      Space below
{\itshape}% Body font
{}%         Indent amount (empty = no indent, \parindent = para indent)
{\bfseries\sffamily}% Thm head font
{.}%        Punctuation after thm head
{ }% Space after thm head: \newline = linebreak
{}%         Thm head spec

\newtheoremstyle{definitionnew}% name
{5pt}%      Space above, empty = `usual value'
{5pt}%      Space below
{\normalfont}% Body font
{}%         Indent amount (empty = no indent, \parindent = para indent)
{\bfseries\sffamily}% Thm head font
{.}%        Punctuation after thm head
{ }% Space after thm head: \newline = linebreak
{}%         Thm head spec

\theoremstyle{plainnew}

\newtheorem{thm}{Theorem}[section]
\newtheorem{prop}[thm]{Proposition}
\newtheorem{lemma}[thm]{Lemma}
\newtheorem{cor}[thm]{Corollary}

\theoremstyle{definitionnew}
\newtheorem{defin}[thm]{Definition}
\newtheorem{esempi}[thm]{Examples}
\newtheorem*{esempi*}{Examples}

\newtheorem*{esempio*}{Example}
\newtheorem{osser}[thm]{Remark}
\newtheorem{noth}[thm]{Notation}

\begin{document}
\newpage
\maketitle
%%%%%%%%%%%%%%%%%%%%%%%%%%%
% abstract, keywords and Subject classification are optional.
%%%%%%%%%%%%%%%%%%%%%%%%%%%
\begin{abstract}
    In this paper, exploiting the work done for my master's thesis, a new construction to associate an internal category to an enriched one is presented. The key concept is that of extensive ambient category and the construction follows the one that associates a category whose idempotents split to a given one. The association turns out to be functorial and left adjoint to an already known one when we restrict to a particular class of internal categories whose idempotents split in some canonical way and impose a size restriction.
\end{abstract}

% Most people don't use these, so they are "commented out"
% by starting the lines with a "%"
%\begin{keywords}
%   \LaTeX, typesetting
%\end{keywords}

%\begin{AMS}
%   50C60, 18C25
%\end{AMS}

%%%%%%%%%%%%%%%%%%%%%%
% % Here is the start of the Text
%%%%%%%%%%%%%%%%%%%%%%
\section{Introduction}
The early motivation for this work comes from \cite{categ} where Brooke-Taylor and Calderoni pose the question of how to \qu{categorify Borel-reducibility} and find the correct categorical stage for this set-theoretical concept. They suggest to look at categories internal to the standard Borel spaces' one, which turns out to be cartesian and countably extensive \cite{chen19}. The concept of countable extensivity is a (countable) generalization of Lawvere's {\itshape extensivity} \cite{LW91, LW92} that was later developed by Carboni, Lack and Walters \cite{ext}. Enriched categories \cite{kelly} and internal ones \cite[][Section B2.3]{elephant} are respectively generalizations of locally small and small categories when we replace \ct{Set} with another suitable ambient category. The question this work deals with is how to connect these two structures, exploiting the infinitary extensivity condition of a cartesian monoidal ambient category \ct{V}, in particular how to generate an internal category once given an enriched one. If $\vcat$ and $\catv$ are respectively the categories of enriched categories and internal ones to $\ct{V}$, a functorial construction $\en\!\,$$\colon {\catv}\to{\vcat}$ is presented by Cottrel, Fujii and Power in \cite{extrel}; there they show also how to build a functor in the other direction that is left adjoint to the first one, exploiting the additional hypothesis of cartesian closedness of $\ct{V}$. This construction has been generalised by Wong \cite{stacks} dropping the cartesian condition of the ambient category and requesting that equalizer preserve coproduct instead of extensivity. In the following pages is presented a functor $\inter$ which associate an enriched category over \ct{V} to one internal to it. The latter is obtained miming the construction of the {\itshape splitting category} or {\itshape idempotent completion}  \cite{catall}, also referred to as the {\itshape Karoubian envelope} \cite[][Section A1.1]{elephant} of a given category. Restricting $\catv$ to the subcategory $\catspv$ of internal categories whose idempotents split in a canonical way and internal functors that preserve it and imposing a size restriction based on \ct{V} to both enriched and internal categories, the functor $\inter$ turns out to be left adjoint to $\en$. This construction heavily relies on the extensivity of the ambient, which is fundamental to express the calculations. The other crucial ingredient is the concept of {\itshape canonical split}: the transposition in the internal setting of two properties some good splits possess; this and the restriction to the internal categories of interest allow to prove the adjunction. 
In section \ref{back} some definitions are recalled together with the basic facts regarding the splitting category; finally the canonical splitting condition is stated and the category $\catspv$ is defined, before concluding the section with a remind of how $\en$ is defined. In section \ref{functor} the functor $\inter$ is constructed and then in section \ref{adj} the main result is presented: the proof of the adjunction $\inter\dashv\en$.

\section{Background}\label{back}
\paragraph{Extensivity} As Cottrell, Fujii and Power \cite{extrel}, Centazzo and Vitale \cite[][Section 4.1]{cenvit03} and Carboni, Lack and Walters \cite{ext} (in a finitary version) do, extensivity will be stated in terms of an equivalence of categories. Let $\ct{V}$ a category, $\kappa$ a cardinal and $(x_i)_{i\in\kappa}$ a $\kappa$-indexed collection of objects of $\ct{V}$. A category with small coproducts will be one with $\kappa$-indexed coproducts for every cardinal $\kappa$. First some other concept that will be useful are defined.

\begin{defin}
	A category with $\kappa$-indexed coproducts and binary products is \dfn{distributive} if the canonical arrow $\coprod_\kappa(Y\times x_i)\to Y\times(\coprod_\kappa x_i)$ is an isomorphism.
\end{defin}

\begin{defin}
	The initial object $0$ is \dfn{strict} if every arrow $x\to 0$ is an isomorphism.
\end{defin}

\begin{prop}
	{\normalfont\bfseries\sffamily(\cite[][Proposition 3.2]{ext})} In a distributive category the initial object is strict.
\end{prop}

\begin{defin}
	In a category with $\kappa$-indexed coproducts they are
	\begin{description}[left=10pt]
		\item[e2\label{e2}]\dfn{universal} if for every arrow $f\colon Y\to\coprod_{i\in\kappa}x_i$, if for every $i\in\kappa$ the following is a pullback square 
		$$\begin{tikzcd}
			y_i\ar[r]\ar[d]&x_i\ar[d]\\ Y\ar[r,"f"]&\coprod_{i\in\kappa}x_i
		\end{tikzcd}$$ then $\coprod_\kappa y_i\to Y$ is an isomorphism;
		\item[e3\label{e3}] \dfn{disjoint} if for every family $(x_i)_{i\in\kappa}$ if $i\neq j$ the following, where $0$ is the initial object, is a pullback square
		$$\begin{tikzcd}
			0\ar[d]\ar[r]&x_i\ar[d]\\x_j\ar[r]&\coprod_\kappa x_i.
		\end{tikzcd}$$
	\end{description}
\end{defin}
If $\ct{V}$ admits $\kappa$-indexed coproducts then there is a functor \begin{equation}\label{1}
	\coprod\colon\prod_{i\in\kappa}(\ct{V}/x_i)\longrightarrow\ct{V}/(\coprod_{i\in\kappa}x_i) 
\end{equation} that sends a collection $(f_i\colon y_i\to x_i)_{i\in\kappa}$ in $\coprod_{i\in\kappa}f_i\colon\coprod_{i\in\kappa}y_i\to\coprod_{i\in\kappa}x_i$.

\begin{osser}
	By definition of coproduct to obtain an arrow $\coprod_{i\in\kappa}y_i\to\coprod_{i\in\kappa}x_i$ we should first compose each $f_i$ with the immersion of $x_i$ in the coproduct. The notation for the induced arrow won't take this into account for semplicity.
\end{osser}

\begin{defin}
	A category $\ct{V}$ that admits $\kappa$-indexed coproduct is $\bm{\kappa}$-\dfn{extensive} if the functor \eqref{1} is an equivalence of categories. 
\end{defin}

\begin{defin}
	A category $\ct{V}$ with small coproducts that for every cardinal $\kappa$ is $\kappa$-extensive is said \dfn{extensive}.
\end{defin}
The following results give two useful characterizations of an extensive category:
\begin{prop}
	{\normalfont\bfseries\sffamily(\cite[][Corollary 2.4 ]{extrel}, \cite[][Section 4.2 Exercise 1]{cenvit03})} A category \ct{V} with small coproducts is extensive if and only if it has all pullbacks along injections and the following two conditions hold:
	\begin{description}[left=10pt] 
		\item[e1\label{e1}] for every small family of arrows $(f_i\colon y_i\to xi)_{i\in\kappa}$ for all $i\in\kappa$ the following square, where the vertical arrows are injections, is a pullback
		$$\begin{tikzcd}
			y_i\ar[r,"f_i"]\ar[d]&x_i\ar[d]\\\coprod_{i\in\kappa}y_i\ar[r,"\coprod_{i\in\kappa}f_i"]&\coprod_{i\in\kappa}x_i;
		\end{tikzcd}$$ 
		\item[e2] coproducts are universal.
	\end{description}
\end{prop}

\begin{cor}
	{\normalfont\bfseries\sffamily(\cite[][Section 4.2 Exercise 2]{cenvit03})} A category \ct{V} with small coproducts is extensive if and only if it has all pullbacks along injections and the following two conditions hold:
	\begin{description}[left=10pt]
		\item[e2] coproducts are universal;
		\item[e3] coproducts are disjoint. 
	\end{description}
\end{cor}
\begin{proof}
	It's straightforward to see that we obtain \ref{e3} from \ref{e1} writing, for $i\neq j$, $f_j=\mathsf{id}_{x_j}$ and for every $k$ different from $j$, $f_k\colon 0 \to x_k$. 
	
	Let's show now that \ref{e1} follows from \ref{e2} and \ref{e3}. Let $z_k$ the pullback of the injection of $x_k$ and the arrow $\coprod_{i\in\kappa}f_i$. By \ref{e2} we have \begin{equation}\label{2}z_k\simeq\coprod_{j\in\kappa}z_k\times_{\coprod y_i}y_j\end{equation} but the injection of $y_j$ followed by $\coprod f_i$ is equal to $f_j$ composed with the injections of $x_j$. Then $z_k\times_{\coprod y_i}y_j$ is isomorphic to the following pullback $$\begin{tikzcd}
		&&x_k\ar[d]\\y_j\ar[r,"f_j"]&x_j\ar[r]&\coprod_{i\in\kappa}x_i
	\end{tikzcd}$$
	which, by \ref{e3} is $y_k$ if $j=k$ and $0$ otherwise. One concludes by \eqref{2}.
\end{proof}

\begin{cor}
	{\normalfont\bfseries\sffamily(\cite[][Section 4.2 Exercise 3]{cenvit03})} In a category \ct{V} with small coproducts and finite limits tfae:
	\begin{description}[left=10pt]
		\item[e2] coproducts are universal;
		\item[e2'\label{e2'}] if, for every $i$ in $\kappa$ the left hand square is a pullback, then the right hand
		square is also a pullback $$\begin{tikzcd}
			y_i\ar[dd]\ar[r]&x_i\ar[d]&[2ex]\coprod_{i\in\kappa}y_i\ar[r]\ar[d]&\coprod_{i\in\kappa}x_i\ar[d]\\&\coprod_{i\in\kappa}x_i\ar[d]&[2ex]A\ar[r]&B\\A\ar[r]&B&[2ex]&
		\end{tikzcd}$$
	\item[e2''\label{e2''}] \begin{enumerate}[i)]\item\label{i} \ct{V} is distributive;
		\item\label{ii} if, for all $i\in\kappa$, $\begin{tikzcd}[cramped, sep=small]e_i\ar[r]& x_i\ar[r,shift left ]\ar[r, shift right]&Y\end{tikzcd}$ is an equalizer then $\begin{tikzcd}[cramped, sep=small]\coprod_\kappa e_i\ar[r]& \coprod_\kappa x_i\ar[r,shift left ]\ar[r, shift right]&Y\end{tikzcd}$ is one as well.
	\end{enumerate}
	\end{description}
\end{cor}
\begin{proof}
	To prove that \ref{e2'} follows from \ref{e2} one obtains an arrow from the vertex of a cone over the cospan we're interested in by expressing it as a coproduct thanks to \ref{e2} and then exploiting the pullback on the left. The vice-versa follows from $A=Y$ and $B=\coprod_{i\in\kappa}x_i$ and the fact that taking pullbacks preserves isomorphisms. Clearly \ref{e2'} implies \ref{i} with $A=Y$ and $B$ the terminal object, since a product is a pullback with the terminal object in the lower right corner, while \ref{ii} follows from \ref{e2} expressing the equalizer of $f,g$ as pullback of $<f,g>$ and $<\mathsf{id}_Y,\mathsf{id}_Y>$. One derives the vice-versa \ref{e2} since as in \cite[][Theorem 1, V.2]{S98} the pullback of the injection of $x_j$ and $f\colon Y\to\coprod_\kappa x_i$ is an equalizer $d_i\to Y\times x_i$ then \ref{e2''} allow to conclude.
\end{proof}

\begin{cor}
	In an extensive category the initial object is strict.
\end{cor}

\begin{cor}\label{inter}
	Let $\ct{V}$ an extensive category, $(x_i)_{i\in\kappa}$ a family of objects, $J_1,J_2$ subsets of $\kappa$. Then the following is a pullback diagram:
	$$\begin{tikzcd}
		\coprod_{i\in J_1\cap J_2}x_i\ar[d]\ar[r]&\coprod_{i\in J_1}x_i\ar[d]\\ \coprod_{i\in J_2}x_i\ar[r]&\coprod_{i\in \kappa}x_i.
	\end{tikzcd}$$
\end{cor}

\begin{proof}
	Apply twice \ref{e2'}: one for the upper-right and one for the lower-left vertex, then the disjointness of coproducts.
\end{proof}

\begin{esempi}
	A big class of extensive categories as examples can be found in \cite{extrel}. The categories $\ct{SBor}$ of standard borel spaces and measurable functions and $\ct{Pol}$ of polish spaces and continuous functions are $\omega$-extensive \cite{chen19}.
\end{esempi}

\paragraph{Splitting category} Let $\ct{C}$ be a category. An \dfn{idempotent} arrow $e$ is such that $ee=e$; domain and codomain of such an arrow will be denoted $\partial e$.

\begin{defin}
	The category \ct{C} \dfn{splits the idempotents} if for every idempotent $e\colon x\to x$ there exist two arrows $r\colon x\to y$ (retraction) and $s\colon y\to x$ (section) such that the following is a commutative diagram
	$$\begin{tikzcd}
		x\ar[dr,"r"]\ar[rr,"e"]&&x\ar[dr,"r"]&\\&y\ar[ur,"s"]\ar[rr,"\mathsf{id}"]&&y.
	\end{tikzcd}$$
\end{defin}

\begin{prop}
	If an idempotent $e$ splits in \ct{C} then the section is an equalizer and the retraction is a coequalizer:
	$$\begin{tikzcd}
		y\ar[r,"s"]&x\ar[r,shift left,"e"]\ar[r,shift right,"\mathsf{id}_x", swap]&x\ar[r,"r"]&y.
	\end{tikzcd}$$
\end{prop}
\begin{proof}
	It's straightforward to verify that the two maps satisfy the respective universal property.
\end{proof}

\begin{osser}
	A split, if it exists, is unique up to (canonical) isomorphism. This will be the reason to ask for a specific split to prove the adjunction. A question will be if this condition is necessary, maybe studying the situation at a $2$-categorical level.
\end{osser}

Let $\idem{\ct{C}}$ the class of idempotent arrow of \ct{C}.

\begin{defin}
	Let\begin{itemize}
		\item $\split{C}_0=\idem{\ct{C}}$;
		\item for $e,d$ in $\split{C}_0$ an arrow $f\colon e\to d$, an arrow $f\colon\partial e\to\partial d$ of \ct{C} such that
		$$\begin{tikzcd}
			\partial e\ar[r,"f"]\ar[d,"e"]&\partial d\\ \partial e\ar[r,"f"] &\partial d\ar[u,"d"]\end{tikzcd}$$
		commutes, equivalently $f=dfe$;
		\item $\split{C}_1$ the class of arrows $f\colon e\to d$ with $e,d$ in $\split{C}_0$.
	\end{itemize}
\end{defin}

\begin{osser}
	If $f\colon e\to d$ is in $\split{C}_1$ then $fe=f$ and $df=f$
\end{osser}

Now some routine computations allow one to check that these data together with the composition inherited from \ct{C}, obvious maps domain and codomain and the choice of the arrow $e\colon e\to e$ as identity for an idempotent $e$ give rise to a category which will be denoted \split{C}.

\begin{prop}\label{split}
	The category \split{\ct{C}} splits the idempotents.
\end{prop}
\begin{proof}
	Let $f\colon e\to e$ an idempotent. It's easy to verify that the following diagram gives a split for $f$
	$$\begin{tikzcd}
		e\ar[dr,"f"]\ar[rr,"f"]&&e\ar[,dr,"f"]&\\&f\ar[ur,"f"]\ar[rr,"f"]&&f.
	\end{tikzcd}$$
\end{proof}

\paragraph{Internal splitting} Let $A=(A_0,A_1,A_2,\partial_1,\partial_2,\iota,c)$ an internal category to \ct{V}, a category with pullbacks. Respectively the notation for an internal category indicates the objects of objects, arrows and composable arrows, domain, codomain, identity and composition; the maps $A_2\to A_1$ arising from the pullback of domain and codomain will be $\Pi_1$ and $\Pi_2$; $A_3$ will be the object of composable triples. We want now to define the notion of splitting of idempotents in the internal setting.

\begin{defin}
	Let \begin{itemize}
		\item (\dfn{endoarrows} object) $(\endo{A},E_A)$ the equalizer of $$\begin{tikzcd}
			A_1\ar[r,shift left,"\partial_0"]\ar[r,shift right, swap, "\partial_1"]&A_0;\end{tikzcd}$$
		\item (\dfn{idempotent} object) $(\idem{A}, I_A)$ the equalizer of
		$$\begin{tikzcd}
			\endo{A}\ar[r,shift left, "c(E_A\times_{A_0}E_A)"]\ar[r,shift right, swap, "\mathsf{id}"]&[6ex]\endo{A}.
		\end{tikzcd}$$
	\end{itemize}
\end{defin}
\begin{defin}
	An internal category $A$ to \ct{V} \dfn{splits the idempotents} if there exist two arrows $R_A,S_A\colon\idem{A}\to A_1$ such that
	\begin{itemize}
		\item $\partial_0R=\partial_0EI=\partial_1S$;
		\item $\partial_1R=\partial_0S$;
		\item $c(S\times_{A_0}R)=EI$;
		\item $c(R\times_{A_0}S)=\iota\partial_0S$.
	\end{itemize}
\end{defin}

\begin{osser}
	The identity arrow $\iota$ induces an arrow $\iota^*\colon A_0\to\idem{A}$.
\end{osser}
\begin{defin}
	An internal category $A$ is said to split \dfn{canonically} the idempotents if it splits the idempotents and
	\begin{itemize}
		\item $R_A\iota^*=S_A\iota^*=\iota$;
		\item for every object $a$ in \ct{V} and arrow $f\times_{A_0}g\times_{A_0}f\colon a \to A_3$ that factors through $EI\times_{A_0}EI\times_{A_0}EI$ and such that the following diagram commutes
		\begin{equation}\label{canon}\begin{tikzcd}
				a\ar[r,"<f^*{,}g^*{,}f^*>"]\ar[dr,"f\times_{A_0}g\times_{A_0}f", swap]&[3em]\idem{A}\times\idem{A}\times\idem{A}\ar[d,"EI\times_{A_0}EI\times_{A_0}EI"]\ar[dr,"EI\pi_2", start anchor={south east}]&[2em]\\&[3em]A_3\ar[r,"c(\mathsf{id}_{A_1}\times_{A_0}c)", swap]&[2em]A_1
		\end{tikzcd}\end{equation}
		one has
		\begin{equation}\label{canon2}
			\begin{array}{r@{\;=\;}l}
				R_A(c(R_Af^*\times_{A_0}c(g\times_{A_0}S_Af^*)))^*&c(R_Ag^*\times_{A_0}S_Af^*)\quad\text{and}\\ S_A(c(R_Af^*\times_{A_0}c(g\times_{A_0}S_Af^*)))^*&c(R_Af^*\times_{A_0}S_Ag^*).\end{array}
		\end{equation}
	\end{itemize}
\end{defin}

\begin{osser}
The last two conditions are respectively the internal analogous of the requests that identities split through themselves and, if we have $fgf=g$ for idempotents $f$ and $g$, then the red arrow 
$$
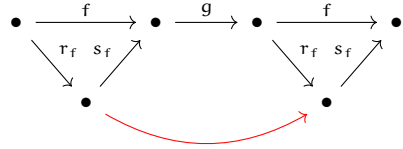
\begin{tikzcd}[column sep=small]
	\bullet\ar[rr,"f"]\ar[dr,"r_f"]&&\bullet\ar[rr,"g"]&&\bullet\ar[rr,"f"]\ar[dr,"r_f"]&&\bullet\\&\bullet\ar[ur,"s_f"]\ar[rrrr,color=red, bend right]&&&&\bullet\ar[ur,"s_f"]&
\end{tikzcd}$$
splits as $r_gs_f$ followed by $r_fs_g$.
\end{osser}

\begin{esempi}
	\begin{itemize}
		\item The split described in Proposition \ref{split} is canonical;
		\item the factorisation of a map in \ct{Set} as the same map with the image as codomain followed by the inclusion of the image in the codomain is a canonical split;
		\item in the full subcategory of \ct{Set} formed by finite cardinals a split such that the section is monotone is canonical.
	\end{itemize}
\end{esempi}

\begin{osser}
	Starting from the split in \ct{Set} just described, for a given identity we define instead the split to be a couple of bijections different from it and for some suitable $f$ and $g$ we define the split of $r_fgs_f$ as $r_gs_f$ followed by a permutation of codomain and the inverse permutation followed by $r_fs_g$. Trivially this split isn't canonical.
\end{osser}

Finally one needs to define an internal functor that preserves the split. 
\begin{osser}
	For an internal functor $F\colon A\to B$ an arrow $F_1^*$ in \ct{V} such that
	$$\begin{tikzcd}
		\idem{A}\ar[d,"F_1^*"]\ar[r,"E_AI_A"]&A_1\ar[d,"F_1"]\\\idem{B}\ar[r,"E_BI_B"]&B_1
	\end{tikzcd}$$
	commutes arises thanks to the commutativity of $F_1$ with compositions and the u.p. of $\idem{B}$.
\end{osser}

\begin{defin}
	An internal functor  $F\colon A\to B$ between internal  categories that split the idempotents \dfn{preserves the split} if 
	$$
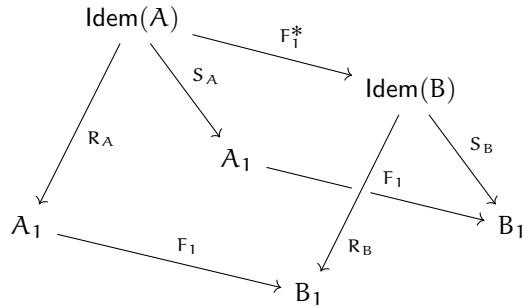
\begin{tikzcd}[row sep=small, column sep=tiny]
		&\idem{A}\ar[drrr,"F_1^*"]\ar[dddl,"R_A"]\ar[ddr,"S_A"]&&&&\\&&&&\idem{B}\ar[ddr,"S_B"]&\\ &&A_1\ar[drrr,"F_1"]&&&\\A_1\ar[drrr,"F_1"]&&&&&B_1\\&& &B_1\ar[from=uuur,"R_B", crossing over, near end]&&
	\end{tikzcd}$$
	commutes.
\end{defin}

It's routine computations to verify that the composition of split preserving functors preserves the split and so the identity.

\begin{defin}\label{catsplit}
	The subcategory of $\catv$ whose object are internal categories that splits idempotents and whose arrows are split preserving internal functors is $\catspv$.
\end{defin}

\paragraph{The functor $\en$}
Given $\ct{V}$ with finite limits let's briefly recall the construction of the $\ct{V}$-enriched category $\en(A)$ for an internal category $A$.
\begin{defin}
	Let
	\begin{itemize}
		\item $\Gamma(A)_0$ the class of arrows $\mathbf{1}\xrightarrow{\alpha}A_0$;
		\item for every $\alpha,\beta\in\Gamma(A)_0$, $(\Gamma(A)(\alpha,\beta),\uparrow_\alpha^\beta,!)$ the limit of the following diagram
		
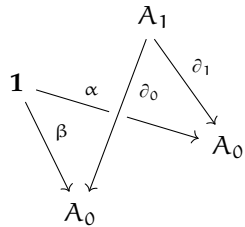
\begin{center}
			\begin{tikzcd}[row sep=small, column sep=tiny]
				&&A_1\ar[ddr,"\partial_1", start anchor={[xshift=-1ex]}]&\\ \mathbf{1}\ar[ddr,"\beta"]\ar[drrr,"\alpha", near start]&&&\\&&&A_0\\ &A_0\ar[from=uuur,"\partial_0", crossing over, near start]&&
			\end{tikzcd}
		\end{center}
	\end{itemize}
\end{defin}

\begin{defin}
	For every $\alpha,\beta,\gamma$ in $\Gamma(A)_0$ define
	\begin{itemize}
		\item $\circ_{\alpha\beta\gamma}$ the unique arrow $\Gamma(A)(\beta,\gamma)\times\Gamma(A)(\alpha,\beta)\to\Gamma(A)(\alpha,\gamma)$ such that $$\uparrow_\alpha^\gamma\circ_{\alpha\beta\gamma}=c(\uparrow_\beta^\gamma \pi_1\times_{A_0}\uparrow_\alpha^\beta \pi_2).$$
		\item $i_\alpha$ the unique arrow induced by $\alpha$ and $\iota\alpha\colon\mathbf{1}\to A_1$, then the arrow $\mathbf{1}\to\Gamma(A)(\alpha,\alpha)$ such that $$\uparrow_\alpha^\alpha i_\alpha=\iota\alpha.$$
	\end{itemize}
\end{defin}
Then it isn't difficult to verify that this data define an enriched category and the construction is funtorial in a natural way.

\section{The functor $\inter$}\label{functor}
From now on \ct{V} will be a cartesian monoidal extensive category with finite limits and small coproducts.

\begin{defin}
	Let $\catvs$ the category of internal categories $A$ to \ct{V} such that the collection of arrows $\mathbf{1}\to A_1$ is a set, $\catspvs$ its subacategory whose objects canonically split the idempotents and whose arrows preserve this; $\vcats$ the category of \ct{V}-enriched categories \ct{M} such that $\ct{M}_0$ is a set and for all $a,b\in\ct{M}_0$ the collection of arrows $\mathbf{1}\to\ct{M}_{ab}$ as well.
\end{defin}

\begin{defin}
	Let $\en'$$\!\colon\catspvs\to\vcats$ the functor defined as $\en$
\end{defin}

\begin{osser}
	All of the following works anyway if \ct{V} has instead only $\kappa$-indexed coproducts, is $\kappa$-extensive and one considers $\catvk$ and $\vcatk$ defined accordingly. 
\end{osser}

\paragraph{The object of objects}
\begin{defin}
	For \ct{M} as above, define $(\inter(\ct{M})_0,\epsilon)$ the equalizer of the following triple
	$$\begin{tikzcd}[column sep=huge]
		\coprod_{a,b,c\in\ct{M}_0}\ct{M}_{bc}\times\ct{M}_{ab}\ar[r,"\mathsf{comp}" description]\ar[r,"\pi_1",shift right=1ex, swap]\ar[r,"\pi_2" ,shift left=1.5ex]&\coprod_{a,b\in\ct{M}_0}\ct{M}_{ab}
	\end{tikzcd}$$
	where $\pi_1,\pi_2$ and $\mathsf{comp}$ are induced respectively by first and second projection and enriched composition.	
\end{defin}

\begin{osser}
	If, for every $a,b\in\ct{M}_0$, we denote $\mathfrak{M}_{ab}$ the following pullback
	$$\begin{tikzcd}
		\mathfrak{M}_{ab}\ar[d]\ar[rr]&&[6ex]\ct{M}_{ab}\ar[d]\\\inter(\ct{M})_0\ar[r,"\epsilon"]&\coprod_{a,b,c\in\ct{M}_0}\ct{M}_{bc}\times\ct{M}_{ab}\ar[r,"\mathsf{comp}" description]\ar[r,"\pi_1",shift right=1ex, swap]\ar[r,"\pi_2" ,shift left=1.5ex]&[6ex]\coprod_{a,b\in\ct{M}_0}\ct{M}_{ab}
	\end{tikzcd}$$
	by universality of coproducts one has $\inter(\ct{M})_0\simeq\coprod_{a,b\in\ct{M}_0}\mathfrak{M}_{ab}.$
\end{osser}

\begin{osser}
	The following is a pullback diagram by \ref{e1} since $A_2^\ct{M}\simeq\coprod_{a,b}\coprod_c\ct{M}_{bc}\times\ct{M}_{ab}$
	$$\begin{tikzcd}
		\coprod_c\ct{M}_{bc}\times\ct{M}_{ab}\ar[r, "\pi_2"]\ar[d]&\ct{M}_{ab}\ar[d]\\\coprod_{a,b,c\in\ct{M}_0}\ct{M}_{bc}\times\ct{M}_{ab}\ar[r,"\pi_2"]&\coprod_{a,b\in\ct{M}_0}\ct{M}_{ab};
	\end{tikzcd}$$
	considering $\pi_1$ and $\mathsf{comp}$ in the upper left corner we get $$\begin{array}{c@{\quad\text{and}\quad}c}
		\coprod_c\ct{M}_{ab}\times\ct{M}_{ca}&\coprod_c\ct{M}_{cb}\times\ct{M}_{ac}.
	\end{array}$$
\end{osser}
\begin{noth}
	For $a,b,c\in\ct{M}_0$, let $\mathfrak{C}_{abc}=\ct{M}_{bc}\times\ct{M}_{ab}$ and $\mathfrak{C}_{ab,}=\coprod_c\ct{M}_{bc}\times\ct{M}_{ab}$.
\end{noth}
\begin{osser}
	The pullback of the following diagram is $\mathfrak{C}_{abc}$ if $b=d$ and $0$ otherwise by Corollary \ref{inter}:
	$$\begin{tikzcd}
		&\mathfrak{C}_{ad,}\ar[d]\\ \mathfrak{C}_{,bc}\ar[r]&\coprod_{a,b,c\in\ct{M}_0}\ct{M}_{bc}\times\ct{M}_{ab}
	\end{tikzcd}$$ 
\end{osser}

\begin{prop}
	Let $\ct{M}$ as above. Then $\inter(\ct{M})_0\simeq\coprod_{a\in\ct{M}_0}\mathfrak{M}_{aa}.$
\end{prop}
\begin{proof}
	Let $a,b\in\ct{M}_0$ and $a\neq b$. All the squares of the following diagram are pullbacks by extensivity and the previous remarks; the dashed arrow is obtained through the universal property of the related pullback:
	$$\begin{tikzcd}[column sep=large, row sep=small]
		&|[yshift=3em]|0\ar[d]&\\&[5ex]|[xshift=2em]|\mathfrak{C}_{ab,}\ar[ddd,color=blue,crossing over]\ar[dr,color=blue,crossing over]&[5ex]\\[-1em]
		\mathfrak{M}_{ab}\ar[uur,dashed]\ar[ru,color=blue,crossing over]\ar[dd,crossing over]\ar[dr,color=green,crossing over]&[5ex]&[5ex]\ct{M}_{ab}\ar[dd, crossing over]\\[-1em]&[5ex]|[xshift=-2em,yshift=3em]|\mathfrak{C}_{,ab}\ar[from=uuu,crossing over]\ar[d,color=green,crossing over]\ar[ru,color=green,crossing over]&[5ex]\\
		\inter(\ct{M})_0\ar[r,"\epsilon"]&[5ex]\coprod_{a,b,c\in\ct{M}_0}\ct{M}_{bc}\times\ct{M}_{ab}\ar[r,"\mathsf{comp}" description,color=red]\ar[r,"\pi_1" description,shift right=1.5ex, color=green]\ar[r,"\pi_2" description,shift left=1.5ex,color=blue]&[5ex]\coprod_{a,b\in\ct{M}_0}\ct{M}_{ab}
	\end{tikzcd}$$
	From the strictness of the initial object one deduces that if $a\neq b$ then $\mathfrak{M}_{ab}\simeq 0$ and the result follows.
\end{proof}

\begin{osser}
	The object $\mathfrak{M}_a:=\mathfrak{M}_{aa}$ is isomorphic to the equalizer of the pair $\begin{tikzcd}[column sep=huge]
			\ct{M}_{aa}\ar[r,"\mathsf{id}", swap, shift right]\ar[r,"\circ_{aaa}<\mathsf{id}{,}\mathsf{id}>", shift left]&\ct{M}_{aa}.
	\end{tikzcd}$ with $\epsilon_a\colon\mathfrak{M}_a\to\ct{M}_{aa}$.
\end{osser}

\paragraph{The object of arrows and composable arrows}
The following definitions will be clearly inspired to the construction of the splitting category; this will be useful to keep in mind not to get lost in the, sometimes cumbersome, notation.
\begin{defin}
	Let
	\begin{itemize}
		\item $\alpha_{ab}\colon\mathfrak{M}_b\times \ct{M}_{ab}\times\mathfrak{M}_a\xrightarrow{(\epsilon_b\times\mathsf{id}\times\epsilon_a)}\ct{M}_{bb}\times \ct{M}_{ab}\times\ct{M}_{aa}\xrightarrow{\mathsf{id}\times\circ_{aab}}\ct{M}_{bb}\times\ct{M}_{ab}\xrightarrow{\circ_{abb}}\ct{M}_{ab};$
		\item $\balpha=\coprod_{a,b}<\pi_1,\alpha_{ab},\pi_3>$ and $(\inter(\ct{M})_1,\cbalpha)$ the equalizer of the pair $$\begin{tikzcd}[column sep= large]
			\coprod_{a,b}\mathfrak{M}_b\times \ct{M}_{ab}\times\mathfrak{M}_a\ar[r,"\balpha", shift left]\ar[r,"\mathsf{id}", shift right, swap]&\coprod_{a,b}\mathfrak{M}_b\times \ct{M}_{ab}\times\mathfrak{M}_a.
		\end{tikzcd}$$
	\end{itemize}
\end{defin}

\begin{osser}
	Clearly this last object is the coproduct of the equalizers of the components and since $\circ_{aaa}<\epsilon_a,\epsilon_a>=\epsilon_a$ it's easy to see that
	$$\begin{array}{r@{}l}\circ_{aab}<\pi_2,\epsilon_a\pi_3>\cbalpha_{ab}=\pi_2\cbalpha_{ab}&\quad\text{and}\\\circ_{abb}<\epsilon_b\pi_1,\pi_2>\cbalpha_{ab}=\pi_2\cbalpha_{ab}&\end{array}$$ in complete analogy to what happens with the arrows in the splitting category.
\end{osser}

\begin{defin}
	Let
	\begin{itemize}
		\item $\iota\colon\inter(\ct{M})_0\to\inter(\ct{M})_1$ the arrow induced by
		$\begin{tikzcd}[column sep=huge]
			\inter(\ct{M})_0\ar[r,"\coprod_a<\mathsf{id}{,}\epsilon_a{,}\mathsf{id}>"]&\coprod_{a,b}\mathfrak{M}_b\times \ct{M}_{ab}\times\mathfrak{M}_a;
		\end{tikzcd}$
		\item $\partial_0\coloneq\coprod_{a,b}\pi_3\cbalpha_{ab}\colon\inter(\ct{M})_1\to\inter(\ct{M})_0$ and $\partial_1\coloneq\coprod_{a,b}\pi_1\cbalpha_{ab}\colon\inter(\ct{M})_1\to\inter(\ct{M})_0;$
		\item $(\inter(\ct{M})_2, \Pi_1,\Pi_2)$ the pullback of $$\begin{tikzcd}
			&\inter(\ct{M})_1\ar[d,"\partial_1"]\\\inter(\ct{M})_1\ar[r,"\partial_0"]&\inter(\ct{M})_0
		\end{tikzcd}$$ and then the coproduct over $a,b,c\in\ct{M}_0$ of objects $\inter(\ct{M})_{abc}$.
	\end{itemize}
\end{defin}

\begin{comment}
\begin{prop}\label{domcod}
	One has $\partial_0\iota=\partial_1\iota=\mathsf{id}_{\inter(\ct{M})_0}$. Equivalently the following diagram commutes: $$\begin{tikzcd}[row sep=tiny]
		&\inter(\ct{M})_1\ar[dl,"\partial_0" , bend right=15, swap]&&\inter(\ct{M})_1\ar[dr,"\partial_1", bend left=15]&\\\inter(\ct{M})_0&&\inter(\ct{M})_0\ar[ll,"\mathsf{id}" description, bend left=30]\ar[rr,"\mathsf{id}" description, bend right=30]\ar[ul,"\iota" ,bend right=15,swap]\ar[ur,"\iota", bend left=15]&&\inter(\ct{M})_0.
	\end{tikzcd}$$
\end{prop}
\begin{proof}
	Just compute on the components.
\end{proof}
\end{comment}

Now the last arrow that remains to define to complete the data for an internal category is composition. 
\begin{osser}
	The two arrows $$\begin{array}{l@{\quad\text{and}\quad}r}
		\cbalpha_{bc}\Pi_1^{abc}\colon \inter(\ct{M})_{abc}\to\mathfrak{M}_c\times\ct{M}_{bc}\times\mathfrak{M}_b& \cbalpha_{ab}\Pi_2^{abc}\colon \inter(\ct{M})_{abc}\to\mathfrak{M}_b\times\ct{M}_{ab}\times\mathfrak{M}_a
	\end{array}$$
	induce $\inter(\ct{M})_{abc}\to(\mathfrak{M}_c\times(\ct{M}_{bc}\times\mathfrak{M}_b))\times(\mathfrak{M}_b\times(\ct{M}_{ab}\times\mathfrak{M}_a))$ which, composing with $$(\mathfrak{M}_c\times(\ct{M}_{bc}\times\mathfrak{M}_b))\times(\mathfrak{M}_b\times(\ct{M}_{ab}\times\mathfrak{M}_a))\xrightarrow{\simeq}\mathfrak{M}_c(((\times\ct{M}_{bc}\times\mathfrak{M}_b)\times(\mathfrak{M}_b\times\ct{M}_{ab}))\times\mathfrak{M}_a)$$ and
	$\mathsf{id}_{c}\times\circ_{abc}<\pi_2,\pi_5>\times\mathsf{id}_a$, gives an arrow $\inter(\ct{M})_{abc}\to\mathfrak{M}_c\times(\ct{M}_{ac}\times\mathfrak{M}_a).$
\end{osser}

\begin{defin}
	In a nutshell, let $$c'_{abc}=(\mathsf{id}_{c}\times\circ_{abc}<\pi_2,\pi_5>\times\mathsf{id}_a)\bigl(<\cbalpha_{bc}\Pi_1^{abc},\cbalpha_{ab}\Pi_2^{abc}>\bigr)$$ omitting the canonical isomorphism between the products.
\end{defin}

\begin{lemma}
	For $a,b,c\in\ct{M}_0$ the arrow $c'_{abc}$ induces $c_{abc}\colon\inter(\ct{M})_{abc}\to\inter(\ct{M})_{ac}$.
\end{lemma}
\begin{proof}
	It's sufficient to check that composition with $\alpha_{ac}$ is the same as taking the second projection. The following diagram commutes by the associativity of the composition of \ct{M}:
	\begin{equation}\label{assscm}
		\begin{tikzcd}
			\inter(\ct{M})_{abc}\ar[ddr,"c'_{abc}", bend right=25, swap, end anchor={[yshift=3pt]}]\ar[r,"<\cbalpha_{bc}\Pi_1{,}\cbalpha_{ab}\Pi_2>"]&[4em](\mathfrak{M}_c\times(\ct{M}_{bc}\times\mathfrak{M}_b))\times(\mathfrak{M}_b\times(\ct{M}_{ab}\times\mathfrak{M}_a))\ar[d,"\simeq"]\ar[r,"\alpha_{ac}\times\alpha_{ab}", shift left]\ar[r, "<\pi_2{,}\pi_5>", shift right, swap]&[2em]\ct{M}_{bc}\times\ct{M}_{ab}\ar[dddl,"\circ_{abc}", bend left]\\&[4em]\mathfrak{M}_c(((\times\ct{M}_{bc}\times\mathfrak{M}_b)\times(\mathfrak{M}_b\times\ct{M}_{ab}))\times\mathfrak{M}_a)\ar[d,"\mathsf{id}\times \circ_{abc}<\pi_2{,}\pi_5>\times\mathsf{id}"]&[2em]\\&[4em]\mathfrak{M}_c\times\ct{M}_{ac}\times\mathfrak{M}_a\ar[d,"\alpha_{ac}"]&[2em]\\&[2em]\ct{M}_{ac}&[2em]
		\end{tikzcd}
	\end{equation} 
	The right hand side is exactly the second projection and one concludes.
\end{proof}

\begin{osser}
	It's just a matter of simple computations and diagrams similar to \eqref{assscm} to conclude that this data yields an internal category \inter(\ct{M}).
\end{osser}

Finally for an enriched functor $\ct{F}\colon\ct{M}\to\ct{N}$ one construct an internal one $\inter(\ct{F})\colon\inter(\ct{M})\to\inter(\ct{N})$ as follows.
\begin{defin}
	Let $\inter(F)$ be the couple formed by
	\begin{itemize}
		\item $\inter(\ct{F})_0$: the coproduct over $a$ in $\ct{M}_0$ of the arrow $\inter(\ct{F})_{0,a}\colon\mathfrak{M}_a\to\mathfrak{N}_{\ct{F}_0a}$ obtained through the following diagram
		$$\begin{tikzcd}[column sep=huge]
			\mathfrak{M}_a\ar[r,"\ct{F}_{aa}\epsilon_a"]&\ct{N}_{\ct{F}_0a\ct{F}_0a}\ar[r,"\mathsf{id}", swap, shift right]\ar[r,"\circ_{\ct{F}_0a\ct{F}_0a\ct{F}_0a}<\mathsf{id}{,}\mathsf{id}>", shift left]&[4em]\ct{N}_{\ct{F}_0a\ct{F}_0a};
		\end{tikzcd}$$
		\item $\inter(\ct{F})_1$: the coproduct over $a,b$ in $\ct{M}_0$ of the arrow $\inter(\ct{F})_{1,ab}\colon\inter(\ct{M})_{ab}\to\inter(\ct{N})_{\ct{F}_0a\ct{F}_0b}$ obtained through the following diagram $$\begin{tikzcd}[column sep= huge]
			\inter(\ct{M})_{ab}\ar[r,"(\inter(\ct{F})_{0,b}\times\ct{F}_{ab}\times\inter(\ct{F})_{0,a})\cbalpha_{ab}"]&[8em]\mathfrak{N}_{\ct{F}_0b}\times \ct{N}_{\ct{F}_0a\ct{F}_0b}\times\mathfrak{N}_{\ct{F}_0a}\ar[r,"\alpha_{\ct{F}_0a\ct{F}_0b}", shift left]\ar[r,"\mathsf{id}", shift right, swap]&[-1em]\mathfrak{N}_{\ct{F}_0b}\times \ct{N}_{\ct{F}_0a\ct{F}_0b}\times\mathfrak{N}_{\ct{F}_0a}.
		\end{tikzcd}$$
	\end{itemize}
\end{defin}

\begin{osser}
	From the properties an enriched functor satisfies is immediate to prove that this couple effectively constitutes an internal functor. 
\end{osser}

The last goal of this section is to show that the internal category associated in this way to an enriched one splits idempotents canonically and every internal functor obtained in this way from an enriched one preserves the split. Again the guide will be the construction of a split for the splitting category.

\begin{osser}
	It is straightforward to note that $\endo{\inter(\ct{M})}\simeq\coprod_{a}\inter(\ct{M})_{aa}$. Then we get $(\idem{\inter(\ct{M})}, I)$ as coproduct over $a$ in $\ct{M}_0$ of the equalizer of 
	$$\begin{tikzcd}[column sep=huge]
		\inter(\ct{M})_{aa}\ar[r,"c_{aaa}(\mathsf{id}\times_{\mathfrak{M}_a}\mathsf{id})",shift left]\ar[r,"\mathsf{id}",shift right,swap]&[3.5em]\inter(\ct{M})_{aa},
	\end{tikzcd}$$
\end{osser}

\begin{noth}
	An arrow (often into an equalizer) into $\mathfrak{M}_a$ or $\idem{\inter(\ct{M})}_a$ obtained from $f_a$ into respectively $\ct{M}_a$ or $\inter(\ct{M})_a$ by universal property will be noted $f^*_a$.
\end{noth}

\begin{osser}
	Beyond the arrow $\iota^*_a\colon\mathfrak{M}_a\to\idem{\inter(\ct{M})}_a$ we have $[\iota^*]_a\colon\mathfrak{M}_a\to\idem{\inter(\ct{M})}_a$ induced by $<i_a!,\epsilon_a,i_a!>^*\colon\mathfrak{M}_a\to\inter(\ct{M})_{aa}$ where $!$ is the unique arrow $\mathfrak{M}_a\to\mathbf{1}$ and $i_a\colon\mathbf{1}\to\ct{M}_{aa}$ represents the identity of $\ct{M}$ as enriched category. Note that also $i_a$ induces $i_a^*\colon\mathbf{1}\to\mathfrak{M}_a$.
\end{osser}

\begin{lemma}
	There exists an arrow $\omega_a\colon\idem{\inter(\ct{M})}_a\to\mathfrak{M}_a$ such that $\pi_2\cbalpha_{aa}I_a\colon\idem{\inter(\ct{M})}_a\to\ct{M}_{aa}$ factors as $\epsilon_a\omega_a$.
\end{lemma}
\begin{proof}
	Let's compute
	$$\begin{array}{r@{\;=\;}l}
		\circ_{aaa}<\mathsf{id},\mathsf{id}>\pi_2\cbalpha_{aa}I_a&\circ_{aaa}<\pi_2\cbalpha_{aa},\pi_2\cbalpha_{aa}>I_a\\[1ex]&\circ_{aaa}<\pi_2\cbalpha_{aa}\Pi^{aaa}_1,\pi_2\cbalpha_{aa}\Pi^{aaa}_2>(\mathsf{id}\times_{a}\mathsf{id})I_a\\[1ex]&\circ_{aaa}<\alpha_{aa}\cbalpha_{aa}\Pi^{aaa}_1,\alpha_{aa}\pi_2\cbalpha_{aa}\Pi^{aaa}_2>(\mathsf{id}\times_{a}\mathsf{id})I_a\\[1ex]&\pi_2c'_{aaa}(\mathsf{id}\times_{a}\mathsf{id})I_a\\[1ex]&\pi_2\cbalpha_{aa}c_{aaa}(\mathsf{id}\times_{a}\mathsf{id})I_a\\[1ex]&\pi_2\cbalpha_{aa}I_a
	\end{array}$$
	and conclude by u.p. of $\mathfrak{M}_a$.
\end{proof}

\begin{defin}
	Let 
	\begin{itemize}
		\item $\Omega=\coprod_{a\in\ct{M}_0}<\pi_1\cbalpha_{aa}I_a,\omega_a,\pi_3\cbalpha_{aa}I_a>$;
		\item $S_{\inter(\ct{M}),a}\colon\idem{\inter(\ct{M})}\to\inter(\ct{M})_{aa}$ the arrow induced by
		\begin{equation}\label{splits}
			\begin{tikzcd}
				\idem{\inter(\ct{M})}_a\ar[r,"\Omega_a"]&\mathfrak{M}_a\times\mathfrak{M}_a\times\mathfrak{M}_a\ar[r,"<\pi_1{,}\epsilon_a\pi_2{,}\pi_2>"]&[4em]\mathfrak{M}_a\times\ct{M}_{aa}\times\mathfrak{M}_a;
			\end{tikzcd}
		\end{equation}
	\item $R_{\inter(\ct{M}),a}\colon\idem{\inter(\ct{M})}\to\inter(\ct{M})_{aa}$ the one induced by
		$$\begin{tikzcd}
			\idem{\inter(\ct{M})}_a\ar[r,"\Omega_a"]&\mathfrak{M}_a\times\mathfrak{M}_a\times\mathfrak{M}_a\ar[r,"<\pi_2{,}\epsilon_a\pi_2{,}\pi_3>"]&[4em]\mathfrak{M}_a\times\ct{M}_{aa}\times\mathfrak{M}_a.
		\end{tikzcd}$$
	\end{itemize}
\end{defin}

\begin{lemma}\label{ristretto}
	One has that $\partial_1R[\iota^*]=\partial_0S[\iota^*]=\mathsf{id}_{\inter(\ct{M})_0}$.
\end{lemma}
\begin{proof}
	On the component, composing with $\epsilon_a$ one has
	\begin{equation}\label{spliteq}\begin{array}{r@{\;=\;}l}
		\epsilon_a\partial_{1,a}R_a[\iota^*]_a&\epsilon_a\pi_1\cbalpha_{aa}R_a[\iota^*]\\[1ex]&\epsilon_a\pi_1<\pi_2{,}\epsilon_a\pi_2{,}\pi_3>\Omega_a[\iota^*]_a\\[1ex]&\epsilon_a\pi_2\Omega_a[\iota^*]_a\\[1ex]&\epsilon_a\omega_a[\iota^*]_a\\[1ex]&\pi_2\cbalpha_{aa}I_a[\iota^*]\\[1ex]&\pi_2\cbalpha_{aa}<i_a!,\epsilon_a,i_a!>^*\\[1ex]&\epsilon_a
	\end{array}\end{equation}
	and we conclude by u.p. of $\mathfrak{M}_a$.
\end{proof}

\begin{prop}
	The couple $(R_\inter(\ct{M}),S_\inter(\ct{M}))$ is a canonical split for $\inter(\ct{M})$.
\end{prop}
\begin{proof}
	Some simple computation yields that the above couple is a split. Let's show briefly that it is canonical.
	From computations similar to \eqref{spliteq} one deduces $R_a\iota^*_a=S_a\iota^*_a=\iota_a$. Now for $f^*,g^*\colon A\to\idem{\inter(\ct{M})}$ as in diagram \eqref{canon}
	\begin{itemize}
		\item let's compute  $R_a(c(Rf^*\times_0c(g\times_0Sf^*)))^*_a$ composing with $\cbalpha_a$:
		\begin{equation}\label{schifus}<\pi_2,\epsilon_a\pi_2,\pi_3>\Omega_a(c_{aaa}(R_af_a^*\times_{\mathfrak{M}_a}c_{aaa}(g_a\times_{\mathfrak{M}_a}S_af_a^*)))^*;\end{equation}
		\item  let's compute $c_{aaa}(Rg^*\times_0Sf^*)_a$ composing with $\cbalpha_a$: $$(\mathsf{id}\times \circ_{aaa}<\pi_2,\pi_5>\times\mathsf{id})<<\pi_2,\epsilon_a\pi_2,\pi_3>\Omega_ag^*,<\pi_1,\epsilon_a\pi_2,\pi_2>\Omega_af^*_a>.$$
	\end{itemize}
	Composing $\pi_1$\eqref{schifus} with $\epsilon_a$ one gets
	$$\epsilon_a\pi_2\Omega_a(c_{aaa}(R_af_a^*\times_{\mathfrak{M}_a}c_{aaa}(g_a\times_{\mathfrak{M}_a}S_af_a^*)))^*=\pi_2\cbalpha_{aa}I_a(c_{aaa}(R_af_a^*\times_{\mathfrak{M}_a}c_{aaa}(g_a\times_{\mathfrak{M}_a}S_af_a^*)))^*$$ that equals to $$\pi_2(\mathsf{id}\times \circ_{aaa}<\pi_2,\pi_5>\times\mathsf{id})<\cbalpha_{aa}\Pi_1,\cbalpha_{aa}\Pi_2>(Rf_a^*\times_0c_{aaa}(g_a\times_0Sf^*_a))$$ which yields 
	$$\begin{array}{r@{\;=\;}l}
		\circ_{aaa}<\pi_2,\pi_5><\cbalpha_{aa}R_af_a^*,\cbalpha_ac_{aaa}(g_a\times_0Sf^*_a)&\circ_{aaa}<\epsilon_a\omega_af_a^*,\circ_{aaa}<\pi_2\cbalpha_{aa}g_a,\epsilon_a\omega_af_a^*>>\\[1ex]&\circ_{aaa}<\pi_2\cbalpha_{aa}f_a,\circ_{aaa}<\pi_2\cbalpha_{aa}g_a,\pi_2\cbalpha_{aa}f_a>>\\[1ex]&\pi_2\cbalpha_{aa}g_a\\[1ex]&\epsilon_a\omega_ag_a^*\\[1ex]&\pi_2\Omega_ag_a^*
	\end{array}$$ 
	by the properties of $f$ and $g$, then the first components coincide.
	Second and third components are proved equals through some similar computations.
\end{proof}

\begin{osser}
	For an enriched functor $\ct{F}$, the internal functor $\inter(\ct{F})$ preserve the (canonical) split and one can quickly verify that $\inter$ gives a functor $\colon\vcats\to\catspvs$.
\end{osser}

\section{The adjunction}\label{adj}
Due to the restriction on the size of enriched and internal categories being considered, for all respectively $\ct{M}$ and $A$ as such, one has that $\catspvs\bigl(\inter(\ct{M}),A\bigr)$ and $\vcats\bigl(\ct{M},\en(A)\bigr)$ are sets. In this final section will be presented a natural bijection between these hom-set.

\paragraph{The map $\Psi$}
Let $F$ an internal functor $\inter(\ct{M})\to A$.
\begin{defin}
	For every $a$ in $\ct{M}_0$ let
	$$\Psi(F)_0(a)=F_{0,a}i_a^*.$$
\end{defin}

\begin{osser}
	For $a,b$ in $\ct{M}_0$ by the definition of $i_a$ and $i_b$ and then the u.p. of $\inter(\ct{M})_{ab}$ there exists a unique dashed arrow that makes the following diagram commute
	$$\begin{tikzcd}
		&\ct{M}_{ab}\ar{d}{\simeq}[swap]{<!{,}\mathsf{id}{,}!>}&&&\\&\mathbf{1}\times\ct{M}_{ab}\times\mathbf{1}\ar[d,dashed]\ar[rr,"i_b^*\times\mathsf{id}\times i_a^*"]\ar[ddl,bend right,"i_b^*\pi_1",swap, near start]&&\mathfrak{M}_b\times\ct{M}_{ab}\times\mathfrak{M}_a\ar[r,"\alpha_{ab}", shift left]\ar[r,"\mathsf{id}", shift right,swap]&\mathfrak{M}_b\times\ct{M}_{ab}\times\mathfrak{M}_a\\&\inter(\ct{M})_{ab}\ar[urr,"\cbalpha_{ab}", swap, near end]\ar[dr,"\partial_{0,ab}", swap]\ar[dl,"\partial_{1,ab}"]\ar[dd,"F_{1,ab}", near end]&&&\\ \mathfrak{M}_b\ar[dd,"F_{0,b}"]&&\mathfrak{M}_a\ar[dd,"F_{0,a}", swap]\ar[from=uul,"i_a^*\pi_3",near end, bend left, crossing over]&&\\ &A_1\ar[dl,"\partial_1"]\ar[dr,"\partial_0", swap]&&&\\A_0&&A_0.&&
	\end{tikzcd}$$
\end{osser}

\begin{defin}
	Let $\nu_{ab}$ be the composition of the dashed arrow and $<!,\mathsf{id},!>$.
\end{defin}

\begin{defin}
	For $a,b$ in $\ct{M}_0$ let $\Psi(F)_{ab}$ the unique arrow $\ct{M}_{ab}\to\en(A)_{\Psi(F)_0a\Psi(F)_o0}$ such that $$\uparrow_{\Psi(F)_0a}^{\Psi(F)_0b}\Psi(F)_{ab}=F_{1,ab}\nu_{ab}.$$
\end{defin}

\begin{thm}
		The map $\Psi\colon\catspvs\bigl(\inter(-),+\bigr)\to\vcats\bigl(-,\en'(+)\bigr)$ is well-defined and natural in its two arguments.
\end{thm}
\begin{proof}
	To prove that it's well-defined one needs to draw the related diagrams expanding them with the arrows mapping limits in their defining diagrams, then exploit the property of an enriched functor after noting that $\omega$ behaves well with composition and identity. From the fact that $\omega$ behaves well wrt $\inter$ one proves naturality.
\end{proof}

\paragraph{The map $\Phi$}
Let $\ct{F}\colon\ct{M}\to\en(A)$ an enriched functor with $A$ that splits canonically the idempotents through $R_A$ and $S_A$. 

For $a\in\ct{M}_0$, the arrow $\idem{\inter(\ct{M})}_a\xrightarrow{\uparrow_{\ct{F}_0a}^{\ct{F}_0a}\ct{F}_{aa}\pi_2\cbalpha_{aa}I_a}A_1$ can be thought as the restriction of the enriched functor to the idempotent of $\inter(\ct{M})$.
\begin{prop}
	The arrow $\uparrow_{\ct{F}_0a}^{\ct{F}_0a}\ct{F}_{aa}\pi_2\cbalpha_{aa}I_a$ factors through an arrow $\ct{F}_{aa}^*\colon\idem{\inter(\ct{M})}_a\to\idem{A}$ such that $$E_AI_A\ct{F}_{aa}^*=\uparrow_{\ct{F}_0a}^{\ct{F}_0a}\ct{F}_{aa}\pi_2\cbalpha_{aa}I_a.$$
\end{prop}
\begin{proof}
	There exist a unique fill for the following triangle
	$$\begin{tikzcd}[column sep=large, row sep= small]
		\idem{\inter(\ct{M})}_a\ar[dr,"\uparrow_{\ct{F}_0a}^{\ct{F}_0a}\ct{F}_{aa}\pi_2\cbalpha_{aa}I_a"]&\\&A_1\\\idem{A}.\ar[ur,"E_AI_A",swap]&
	\end{tikzcd}$$
	Indeed it's clear that there exist $e_a\colon\ct{M}_{aa}\to\endo{A}$ such that
	$$\begin{tikzcd}[column sep=large, row sep= small]
		\ct{M}_{aa}\ar[dd,"e_a"]\ar[dr,"\uparrow_{\ct{F}_0a}^{\ct{F}_0a}\ct{F}_{aa}"]&\\&A_1\\\endo{A}\ar[ur,"E_A",swap]&
	\end{tikzcd}$$
	commutes; then one proves 
	$c(E_A\times_{A_0}E_A)e_a\pi_2\cbalpha_{aa}I_a=\uparrow_{\ct{F}_0a}^{\ct{F}_0a}\ct{F}_{aa}\pi_2\cbalpha_{aa}I_a$ exploiting the definition of $\idem{\inter(\ct{M})}$ and the fact that composition and $\uparrow$ commute.
\end{proof}

\begin{defin}
	Let $\Psi(\ct{F})_{0,a}\colon\inter(\ct{M})_a\to A_0$ the arrow $\partial_1R_A\ct{F}_{aa}^*\iota_a^*.$
\end{defin}

\begin{lemma}
	For $a,b\in\ct{M}$, given
	\begin{center}
		%\tikzsetnextfilename{frecciafrecce}
		\begin{tikzpicture}[p/.style={rectangle, inner sep=1pt}]
			\node[p] (mcab) at (0,0) {$\inter(\ct{M})_{ab}$};
			
			\node[p] (mab1) at ($(mcab.east)+(3cm,0)$) {$\ct{M}_{ab}$}
			edge[<-, in=0, out=180] node[auto] {$\cbalpha_{ab}$} (mcab.east);
			\node[p] (times1) at ($(mab1.north)+(0,2mm)$) {$\times$};
			\node[p] (ma1) at ($(times1.north)+(0,2mm)$) {$\mathfrak{M}_a$};
			\node[p] (times2) at ($(mab1.south)+(0,-2mm)$) {$\times$};
			\node[p] (mb1) at ($(times2.south)+(0,-2mm)$) {$\mathfrak{M}_b$};
			
			\node[p] (mab2) at ($(mab1.east)+(3cm,0)$) {$\ct{M}_{ab}$}
			edge[<-, in=0, out=180] node[rectangle, fill=white] {$\pi_2$} (mab1.east);
			\node[p] (ma2) at ($(ma1.east)+(3cm,0)$) {$\mathfrak{M}_a$}
			edge[<-] node[above] {$\pi_3$} (mab1.east);
			\node[p] (mb2) at ($(mb1.east)+(3cm,0)$) {$\mathfrak{M}_b$}
			edge[<-] node[below] {$\pi_1$} (mab1.east);
			
			\node[p] (a1) at ($(ma2.east)+(3cm,0)$) {$A_1$}
			edge[<-, in=0, out=180] node[above] {$S_A\ct{F}_{aa}^*\iota^*_a$} (ma2.east);
			\node[p] (a2) at ($(mab2.east)+(2.9cm,0)$) {$A_1$}
			edge[<-, in=0, out=180] node[rectangle, fill=white] {$\uparrow_{\ct{F}_0a}^{\ct{F}_0b}\ct{F}_{ab}$} (mab2.east);
			\node[p] (a3) at ($(mb2.east)+(3cm,0)$) {$A_1$}
			edge[<-, in=0, out=180] node[auto] {$R_A\ct{F}_{bb}^*\iota^*_b$} (mb2.east);
		\end{tikzpicture}
	\end{center}
	the first two arrows are composable and so are the second two.
\end{lemma}
\begin{proof}
	Clear from the definition of $\uparrow$ and the fact that $\partial_1S_A=\partial_1E_AI_A$ and $\partial_0R_A=\partial_0E_AI_A$.
\end{proof}

\begin{cor}
	There exists a unique filling $\Delta^\ct{F}_{ab}\colon\inter(\ct{M})\to A_3$ for
	$$\begin{tikzcd}[column sep=large]
		\inter(\ct{M})_{ab}\ar[drr,"\uparrow_{\ct{F}_0a}^{\ct{F}_0b}\ct{F}_{ab}\pi_2\cbalpha_{ab}\times_{A_0}S_A\ct{F}_{aa}^*\iota^*_a\pi_3\cbalpha_{ab}", bend left]\ar[ddr,"R_A\ct{F}_{bb}^*\iota^*_b\pi_1\cbalpha_{ab}\times_{A_0}\uparrow_{\ct{F}_0a}^{\ct{F}_0b}\ct{F}_{ab}\pi_2\cbalpha_{ab}", bend right, swap]&&\\&A_3\pullback  \ar[d,"<\Pi_1'{,}\Pi_2'>", swap]\ar[r,"<\Pi_2'{,}\Pi_3'>"]&A_2\ar[d,"\Pi_1"]\\&A_2\ar[r, "\Pi_2"]&A_1
	\end{tikzcd}$$
\end{cor}

\begin{defin}
	Let $\Phi(\ct{F})_1$ the coproduct over $a$ in $\ct{M}_0$ of the arrow $\Phi(\ct{F})_{1,ab}\colon\inter(\ct{M})\to A_1$ defined as $c(\mathsf{id}_{A_1}\times_{A_0}c)\Delta^\ct{F}_{ab}.$
\end{defin}

\begin{prop}
	The pair $(\Phi(\ct{F})_0,\Phi(\ct{F})_1)$ defines an internal functor $\Phi(\ct{F})\colon\inter(\ct{M})\to A$.
\end{prop}
\begin{proof}
	Commutativiy with domain and codomain is a matter of simple computations. Moving forward at composition one has
	$$c(\Phi(\ct{F})_{1,bc}\Pi_1^{abc}\times_{A_0}\Phi(\ct{F})_{1,ab}\Pi_2^{abc})O=c(c(\mathsf{id}_{A_1}\times_{A_0}c)\Pi_1^{abc}\Delta^\ct{F}_{bc}\times_{A_0}c(\mathsf{id}_{A_1}\times_{A_0}c)\Pi_2^{abc}\Delta^\ct{F}_{ab})$$ which is the composition of six arrows. By the associativity of composition of $A$ and the fact that $c(S_A\times_{A_0}R_A)=E_AI_A$, the four central amount to $$c(\uparrow_{\ct{F}_0b}^{\ct{F}_0c}\ct{F}_{bc}\pi_2\cbalpha_{bc}\Pi_1^{abc}\times_{A_0}c(E_AI_A\ct{F}_{bb}^*\iota^*_b\pi_1\cbalpha_{ab}\Pi_2^{abc}\times_{A_0}\uparrow_{\ct{F}_0a}^{\ct{F}_0b}\ct{F}_{ab}\pi_2\cbalpha_{ab}\Pi_2^{abc})),$$ which by definition of $\ct{F}_{bb}^*$ equals $$c(\uparrow_{\ct{F}_0b}^{\ct{F}_0c}\ct{F}_{bc}\pi_2\cbalpha_{bc}\Pi_1^{abc}\times_{A_0}c(\uparrow_{\ct{F}_0b}^{\ct{F}_0b}\ct{F}_{bb}\pi_2\cbalpha_{bb}I_b\iota^*_b\pi_1\cbalpha_{ab}\Pi_2^{abc}\times_{A_0}\uparrow_{\ct{F}_0a}^{\ct{F}_0b}\ct{F}_{ab}\pi_2\cbalpha_{ab}\Pi_2^{abc})),$$ and since $\pi_2\cbalpha_{bb}I_b\iota^*_b=\epsilon_b$ one has $$c(\uparrow_{\ct{F}_0b}^{\ct{F}_0c}\ct{F}_{bc}\pi_2\cbalpha_{bc}\Pi_1^{abc}\times_{A_0}c(\uparrow_{\ct{F}_0b}^{\ct{F}_0b}\ct{F}_{bb}\epsilon_b\pi_1\cbalpha_{ab}\Pi_2^{abc}\times_{A_0}\uparrow_{\ct{F}_0a}^{\ct{F}_0b}\ct{F}_{ab}\pi_2\cbalpha_{ab}\Pi_2^{abc})),$$ from which, by definition of composition in $\en(A)$ and since $\ct{F}$ is an enriched functor one gets 
	$$\begin{array}{r@{\;=\;}l}
		c(\uparrow_{\ct{F}_0b}^{\ct{F}_0c}\ct{F}_{bc}\pi_2\cbalpha_{bc}\Pi_1^{abc}\times_{A_0}\uparrow_{\ct{F}_0a}^{\ct{F}_0b}\ct{F}_{ab}\pi_2\cbalpha_{ab}\Pi_2^{abc})&\uparrow_{\ct{F}_0a}^{\ct{F}_0c}\ct{F}_{ac}\circ_{abc}<\pi_2,\pi_5><\cbalpha_{bc}\Pi_1^{abc},\cbalpha_{ab}\Pi_2^{abc}>\\[1ex]&\uparrow_{\ct{F}_0a}^{\ct{F}_0c}\ct{F}_{ac}\pi_2c'_{abc}\\[1ex]&\uparrow_{\ct{F}_0a}^{\ct{F}_0c}\ct{F}_{ac}\pi_2\cbalpha_{ac} c_{abc}.
	\end{array}$$ 
	The first factor being
	$$R_A\ct{F}_{cc}^*\iota^*_c\pi_1\cbalpha_{bc}\Pi_1^{abc}=R_A\ct{F}_{cc}^*\iota^*_c\pi_1\cbalpha_{bc}c_{abc}$$ and the last being $$S_A\ct{F}_{aa}^*\iota^*_a\pi_3\cbalpha_{ab}\Pi_2=S_A\ct{F}_{aa}^*\iota^*_a\pi_3\cbalpha_{ab}c_{abc}$$
	by the defining property of $\inter(\ct{M})_2$, one concludes.
	Finally identity follows again from the associativity of composition in $A$, the definition of $\ct{F}_{bb}^*$ and the properties of the split maps.
\end{proof}

\begin{lemma}
	 If $A$ splits canonically the idempotents the for every enriched functor $\ct{F}\colon\ct{M}\to\en'(A)$, $\Phi(\ct{F})$ preserves the split.
\end{lemma}
\begin{proof}
	Let's compute the composition of $\Phi(\ct{F})_1$ and $S_{\inter(\ct{M})}$:
	$$\begin{array}{r@{\;=\;}l}
		(\Phi(\ct{F})_1S_{\inter(\ct{M})})_a&c(\mathsf{id}_{A_1}\times_{A_0}c)\Delta_{aa}^\ct{F}S_{\inter(\ct{M}),a}\\[1ex]&c(R_A\ct{F}_{bb}^*\iota^*_a\pi_1\cbalpha_{aa}\times_{A_0}c(\uparrow_{\ct{F}_0a}^{\ct{F}_0a}\ct{F}_{aa}\pi_2\cbalpha_{aa}\times_{A_0}S_A\ct{F}_{aa}^*\iota^*_a\pi_3\cbalpha_{aa}))S_{\inter(\ct{M}),a}\\[1ex]&c(R_A\ct{F}_{aa}^*\iota^*_a\pi_1\Omega_a\times_{A_0}c(\uparrow_{\ct{F}_0a}^{\ct{F}_0a}\ct{F}_{aa}\epsilon_a\pi_2\Omega_a\times_{A_0}S_A\ct{F}_{aa}^*\iota^*_a\pi_2\Omega_a))\\[1ex]&c(R_A\ct{F}_{aa}^*\iota^*_a\pi_1\cbalpha_{aa}I_a\times_{A_0}c(\uparrow_{\ct{F}_0a}^{\ct{F}_0a}\ct{F}_{aa}\pi_2\cbalpha_{aa}I_a\times_{A_0}S_A\ct{F}_{aa}^*\iota^*_a\pi_2\Omega_a))\end{array}$$
	by definition of $\Omega_a$ and $S_{\inter(\ct{M}),a}$. Moreover since
	$$\begin{array}{r@{\;=\;}l}
		E_AI_A\ct{F}_{aa}^*\iota^*_a\pi_2\Omega_a&\uparrow_{\ct{F}_0a}^{\ct{F}_0a}\ct{F}_{aa}\pi_2\cbalpha_{aa}I_a\iota^*_a\pi_2\Omega_a\\[1ex]&\uparrow_{\ct{F}_0a}^{\ct{F}_0a}\ct{F}_{aa}\pi_2\cbalpha_{aa}I_a\\[1ex]&E_AI_A\ct{F}_{aa}^*
	\end{array}$$
	one gets $\ct{F}_{aa}^*\iota^*_a\pi_2\Omega_a=\ct{F}_{aa}^*$ by the universal property, and then
	$$\begin{array}{r@{\;=\;}l}
		(\Phi(\ct{F})_1S_{\inter(\ct{M})})_a&c(R_A\ct{F}_{aa}^*\iota^*_a\pi_1\cbalpha_{aa}I_a\times_{A_0}c(\uparrow_{\ct{F}_0a}^{\ct{F}_0a}\ct{F}_{aa}\pi_2\cbalpha_{aa}I_a\times_{A_0}S_A\ct{F}_{aa}^*))\\[1ex]&c(R_A\ct{F}_{aa}^*\iota^*_a\pi_1\cbalpha_{aa}I_a\times_{A_0}c(E_AI_A\ct{F}_{aa}^*\times_{A_0}S_A\ct{F}_{aa}^*))\\[1ex]&c(R_A\ct{F}_{aa}^*\iota^*_a\pi_1\cbalpha_{aa}I_a\times_{A_0}c(c(S_A\times_{A_0}R_A)\ct{F}_{aa}^*\times_{A_0}S_A\ct{F}_{aa}^*))\\[1ex]&c(R_A\ct{F}_{aa}^*\iota^*_a\pi_1\cbalpha_{aa}I_a\times_{A_0}c(S_A\times_{A_0}\iota\partial_1S_A)\ct{F}_{aa}^*)\\[1ex]&c(R_A\ct{F}_{aa}^*\iota^*_a\pi_1\cbalpha_{aa}I_a\times_{A_0}S_A\ct{F}_{aa}^*)
	\end{array}$$
	by the properties of a split and those defining an internal category. Similar computations yields
	$$(\Phi(\ct{F})_1R_{\inter(\ct{M})})_a=c(R_A\ct{F}_{aa}^*\times_{A_0}S_A\ct{F}_{aa}^*\iota^*_a\pi_3\cbalpha_{aa}I_a).$$
	Consider now $$\begin{array}{cl}
		f^*:=\ct{F}_{aa}^*\iota^*_a\pi_1\cbalpha_{aa}I_a&\\g^*:=\ct{F}^*_{aa}&\colon\idem{\inter(\ct{M})}\to\idem{A}\\h^*:=\ct{F}_{aa}^*\iota^*_a\pi_3\cbalpha_{aa}I_a&\end{array}$$
	 and note that $\iota^*_a\pi_1\cbalpha_{aa}I_a=\iota^*_a\pi_3\cbalpha_{aa}I_a$ then $f^*=h^*$. Since $A$ splits canonically the idempotent, if the composition of $E_AI_A f^*, E_AI_A g^*$ and $E_AI_A h^*$ gives the middle one then one concludes, since for $f$ and $g$ respectively $E_AI_A f^*$ and $E_AI_A g^*$, $R_A\Phi(\ct{F})^*_a$ is exactly the LHS and $(\Phi(\ct{F})_1R_{\inter(\ct{M})})_a$ the RHS of the first equation in \eqref{canon2}. The related composition finally is
	$$c(E_AI_A\ct{F}_{aa}^*\iota^*_a\pi_1\cbalpha_{aa}I_a\times_0c(E_AI_A\ct{F}^*_{aa}\times_0E_AI_A\ct{F}_{aa}^*\iota^*_a\pi_3\cbalpha_{aa}I_a))$$ which equals to 
	$$\uparrow_{\ct{F}_0a}^{\ct{F}_0a}\ct{F}_{aa}\circ_{aaa}<\epsilon_a\pi_1\cbalpha_{aa}I_a,\circ_{aaa}<\pi_2\cbalpha_{aa}I_a,\epsilon_a\pi_3\cbalpha_{aa}I_a>>=\uparrow_{\ct{F}_0a}^{\ct{F}_0a}\ct{F}_{aa}\pi_2\cbalpha_{aa}I_a=E_AI_A\ct{F}^*_{aa}.$$
\end{proof}

\begin{cor}
	For every $\ct{V}$-category $\ct{M}$ e internal category $A$ which canonically splits idempotents, the association $\ct{F}\mapsto\Phi(\ct{F})$, given an internal functor $\ct{F}\colon\ct{M}\to\en'(A)$, defines a map $$\Phi\colon\ct{C}\text{-}\mathsf{cat}(\ct{M},\en'(A))\to\mathsf{cat}(\ct{C})_\ct{Sp}(\inter(\ct{M}),A).$$
\end{cor}

\begin{osser}\label{iotastella}
	For an enriched functor $\ct{F}\colon\ct{M}\to\en'(A)$ and $a\in\ct{M}_0$ observe that one has $\ct{F}^*\iota^*_ai_a^*=\iota^*\ct{F}_0a$ since compositions with $E_AI_A$ coincide.
\end{osser}

\begin{thm}
	For every each enriched functor $\ct{F}\colon\ct{M}\to\en'(A)$ one has that $\Psi(\Phi(\ct{F}))=\ct{F}$.
\end{thm}
\begin{proof}
	Compute the action on the set of objects, for $a\in\ct{M}_0$ one has:
	$$\begin{array}{r@{\;=\;}l}
		\Psi(\Phi(\ct{F}))_0a&\Phi(\ct{F})_{0,a}i_a^*\\[1ex]&\partial_1R_A\ct{F}_{aa}^*\iota^*_ai_a^*\\[1ex]&\partial_1R_A\iota^*\ct{F}_0a\\[1ex]&\ct{F}_0a
	\end{array}$$ by remark \ref{iotastella} and because $A$ canonically splits idempotents. For $a,b\in\ct{M}_0$ one gets
	$$\begin{array}{r@{\;=\;}l}
		\uparrow_{\ct{F}_0a}^{\ct{F}_0b}\Psi(\Phi(\ct{F}))_{ab}&\Phi(\ct{F})_{1,ab}\nu_{ab}\\[1ex]& c(\mathsf{id}_{A_1}\times_{A_0}c)\Delta^\ct{F}_{ab}\nu_{ab}\\[1ex]&c(R_A\ct{F}_{bb}^*\iota^*_b\pi_1\cbalpha_{ab}\times_0c(\uparrow_{\ct{F}_0a}^{\ct{F}_0b}\ct{F}_{ab}\pi_2\cbalpha_{ab}\times_0S_A\ct{F}_{aa}^*\iota^*_a\pi_3\cbalpha_{ab}))\nu_{ab}\\[1ex]&c(R_A\ct{F}_{bb}^*\iota^*_bi_b^*!\times_0c(\uparrow_{\ct{F}_0a}^{\ct{F}_0b}\ct{F}_{ab}\times_0S_A\ct{F}_{aa}^*\iota^*_ai_a^*!))\\[1ex]&c(\iota\ct{F}_0b!\times_0c(\uparrow_{\ct{F}_0a}^{\ct{F}_0b}\ct{F}_{ab}\times_0\iota\ct{F}_0a!))\\[1ex]&c(\uparrow_{\ct{F}_0b}^{\ct{F}_0b}i_{\ct{F}_0b}!\times_0c(\uparrow_{\ct{F}_0a}^{\ct{F}_0b}\ct{F}_{ab}\times_0\uparrow_{\ct{F}_0a}^{\ct{F}_0a}i_{\ct{F}_0a}!))\\[1ex]&\uparrow_{\ct{F}_0a}^{\ct{F}_0b}\ct{F}_{ab}
	\end{array}$$ by the remark above and the definintion of identity maps for $\en'(A)$, this let one conclude.
\end{proof}

\begin{cor}
	The map $\Psi$ is surjective.
\end{cor}

\begin{osser}
	One has that $\nu_{aa}\epsilon_a=<i_a!,\epsilon_a,i_a!>^*$ by composing with $\cbalpha_{aa}^*$.
\end{osser}

\begin{lemma}\label{cristretto}
	If $F\colon\inter(\ct{M})\to A$ is an internal functor then 
	$$\begin{tikzcd}
		\mathfrak{M}_a\ar[d, "\iota^*_a"]\ar[r, "{[}\iota_a^*{]}"]&\idem{\inter(\ct{M})}\ar[d,"F_{1,aa}^*"]\\\idem{\inter(\ct{M})}\ar[r, "(\Psi(F)_{aa})^*"]&\idem{A}
	\end{tikzcd}$$
	commutes, meaning $$(\Psi(F)_{aa})^*\iota^*_a=F_{1,aa}^*[\iota^*_a].$$
\end{lemma}
\begin{proof}
	Observe that 
	$$\begin{array}{r@{\;=\;}l}
		E_AI_A(\Psi(F)_{aa})^*\iota^*_a&\uparrow_{\Psi(F)_0a}^{\Psi(F)_0a}\Psi(F)_{aa}\pi_2\cbalpha_{aa}I_a\iota^*_a\\[1ex]&F_{1,aa}\nu_{aa}\pi_2\cbalpha_{aa}<\mathsf{id},\epsilon_a,\mathsf{id}>^*\\[1ex]&F_{1,aa}\nu_{aa}\pi_2<\mathsf{id},\epsilon_a,\mathsf{id}>\\[1ex]&F_{1,aa}\nu_{aa}\epsilon_a\\[1ex]&F_{1,aa}<i_a!,\epsilon_a,i_a!>^*\\[1ex]&F_{1,aa}I_a[\iota^*_a]\\[1ex]&E_AI_AF_{1,aa}^*[\iota^*_a].
	\end{array}$$
	by the above remark and definition of $F_1^*$.
\end{proof}

\begin{lemma}\label{lastbutnotleast}
	If $F\colon\inter(\ct{M})\to A$ is an internal functor then \begin{itemize}
		\item $c_{aab}\bigl(\nu_{ab}\pi_2\times_{\inter(\ct{M})_a}S_{\inter(\ct{M}),a}[\iota^*_a]\pi_3)\bigr)\cbalpha_{ab}=(<i_b!,\pi_2,\pi_3>\cbalpha_{ab}^*)^*$ and
		\item $c_{abb}\bigl(R_{\inter(\ct{M}),b}[\iota^*_b]\pi_1\cbalpha_{ab}\times_{\inter(\ct{M})_b}(<i_b!,\pi_2,\pi_3>\cbalpha_{ab})^*\bigr)=\mathsf{id}_{\inter(\ct{M})_{ab}}.$
	\end{itemize}
\end{lemma}
\begin{proof}
	For the first one compose with $\cbalpha_{ab}$:
	$$\begin{array}{r@{\;=\;}l}
		\cbalpha_{ab}c_{aab}\bigl(\omega_{ab}\pi_2\times_{\inter(\ct{M})_0}S_{\inter(\ct{M})}[\iota^*_a]\pi_3)\cbalpha_{ab}&c'_{aab}\bigl(\nu_{ab}\pi_2\times_{\inter(\ct{M})_0}S_{\inter(\ct{M}),a}[\iota^*_a]\pi_3\bigr)\cbalpha_{ab}\\[1ex]&(\mathsf{id}_b\times \circ_{aab}<\pi_2,\pi_5>\times\mathsf{id}_a)\bigl(<\cbalpha_{ab}\nu_{ab}\pi_2,\cbalpha_{aa}S_{\inter(\ct{M}),a}[\iota^*_a]\pi_3>\bigr)\cbalpha_{ab}\\[1ex]&<i_b!,\pi_2,\pi_3>\cbalpha_{ab}
	\end{array}$$
	since by definition of $\omega_{ab}$ and $S_{\inter(\ct{M}),a}$
	\begin{itemize}
		\item the first component is $\mathsf{id}_bi_b!\pi_2\cbalpha_{ab}=i_b!$;
		\item the middle one is $$\begin{array}{r@{\;=\;}l}\circ_{aab}<\pi_2\cbalpha_{ab}\nu_{ab}\pi_2,\pi_2\cbalpha_{aa}S_{\inter(\ct{M}),a}[\iota^*_a]\pi_3>&\circ_{aab}<\pi_2,\epsilon_a\omega_a[\iota_a^*]\pi_3>\cbalpha_{ab}\\[1ex]&\circ_{aab}<\pi_2,\pi_2\cbalpha_{aa}I_a[\iota_a^*]\pi_3>\cbalpha_{ab}\\[1ex]&\circ_{aab}<\pi_2,\epsilon_a\pi_3>\cbalpha_{ab}\\[1ex]&\pi_2\cbalpha_{ab};\end{array}$$
		\item the third one is $\omega_a[\iota_a]\pi_3\cbalpha_{ab}=\pi_3\cbalpha_{ab}$ by definition of $S$ and analogous computations.
	\end{itemize}
	The second one follows the same way.
\end{proof}

\begin{thm}
	If $F\colon\inter(\ct{M})\to A$ is an internal functor then $\Phi(\Psi(F))=F$. 
\end{thm}
\begin{proof}
	The component $a$ of $\Phi(\Psi(F))_0$ is
	$$\begin{array}{r@{\;=\;}l}
		\Phi(\Psi(F))_{0,a}&\partial_1R_A(\Psi(F)_{aa})^*\iota^*_a\\[1ex]&\partial_1R_AF_{1,aa}^*[\iota^*_a]\\[1ex]&\partial_1F_{1,aa}R_{\inter(\ct{M}),a}[\iota^*_a]\\[1ex]&F_{0,a}
	\end{array}$$
	by Lemma \ref{cristretto} and since $\partial_{1,a}R_{\inter(\ct{M}),a}[\iota^*]_a=\mathsf{id}_{\inter(\ct{M})_a}$ by Lemma \ref{ristretto}.\\
	Component $ab$ of $\Phi(\Psi(F))_1$ is
	$$\begin{array}{r@{\;=\;}l}
		\Phi(\Psi(F))_{1,ab}&c(\mathsf{id}_{A_1}\times_{A_0}c)\Delta^{\Psi(F)}_{ab}\\[1ex]&c\bigl(R_{A,\Psi(F)0b}(\Psi(F)_{bb})^*\iota^*_b\pi_1\times_{A_0}c(\uparrow_{\Psi(F)_0a}^{\Psi(F)_0b}\Psi(F)_{ab}\pi_2\times_{A_0}S_A(\Psi(F)_{aa})^*\iota^*_a\pi_3)\bigr)\cbalpha_{ab}\\[1ex]&c\bigl(R_AF_{1,bb}^*[\iota^*_b]\pi_1\times_{A_0}c(F_{1,ab}\nu_{ab}\pi_2\times_{A_0}S_AF_{1,aa}^*[\iota^*_a]\pi_3)\bigr)\cbalpha_{ab}\\[1ex]&c\bigl(F_{1,bb}R_{\inter(\ct{M}),b}[\iota^*_b]\pi_1\times_{A_0}F_{1,ab}c_{aab}(\nu_{ab}\pi_2\times_{\inter(\ct{M})_a}S_{\inter(\ct{M}),a}[\iota^*_a]\pi_3)\bigr)\cbalpha_{ab}\\[1ex]&c\bigl(F_{1,bb}R_{\inter(\ct{M}),b}[\iota^*_b]\pi_1\cbalpha_{ab}\times_{A_0}F_{1,ab}(<i_b!,\pi_2,\pi_3>\cbalpha_{ab})^*\bigr)\\[1ex]&F_{1,ab}c_{abb}\bigl(R_{\inter(\ct{M}),b}[\iota^*_b]\pi_1\cbalpha_{ab}\times_{\inter(\ct{M})_b}(<i_b!,\pi_2,\pi_3>\cbalpha_{ab})^*\bigr)\\[1ex]&F_{1,ab}
	\end{array}$$
	again by Lemmas \ref{cristretto} and \ref{lastbutnotleast}, and so one concludes.
\end{proof}

\begin{cor}
	The map $\Psi$ is injective, then a natural bijection.
\end{cor}

\begin{thm}
	There is an adjunction $\inter\dashv\en'$.
\end{thm}

\printbibliography[heading=bibintoc]

\end{document}